\documentstyle[12pt]{article}

\newcommand{\arrow}{\rightarrow}
\newcommand{\bairespace}{\omega^{\omega}}
\newcommand{\cantorspace}{2^{\omega}}
\newcommand{\cross}{\times}
\newcommand{\Intersect}{\bigcap}
\newcommand{\intersect}{\cap}
\newcommand{\proof}{\par\noindent proof:\par}
\newcommand{\rationals}{{\Bbb Q}}
\newcommand{\reals}{{\Bbb R}}

\newcommand{\union}{\cup}

\newtheorem{theorem}{Theorem}
\newtheorem{lemma}[theorem]{Lemma}

\newcommand{\res}{\upharpoonright}
\newcommand{\qed}{\nopagebreak\par\noindent\nopagebreak$\blacksquare$\par}

\catcode`\@=11

\font\teneuf=eufm10  scaled 1200
\font\seveneuf=eufm7 scaled 1200
\font\fiveeuf=eufm5  scaled 1200
\font\tenmsx=msxm10  scaled 1200
\font\sevenmsx=msxm7 scaled 1200
\font\fivemsx=msxm5  scaled 1200
\font\tenmsy=msym10  scaled 1200
\font\sevenmsy=msym7 scaled 1200
\font\fivemsy=msym5  scaled 1200

\newfam\euffam
\newfam\msxfam
\newfam\msyfam

\textfont\euffam=\teneuf \scriptfont\euffam=\seveneuf
 \scriptscriptfont\euffam=\fiveeuf
\textfont\msxfam=\tenmsx  \scriptfont\msxfam=\sevenmsx
  \scriptscriptfont\msxfam=\fivemsx
\textfont\msyfam=\tenmsy  \scriptfont\msyfam=\sevenmsy
  \scriptscriptfont\msyfam=\fivemsy

\def\frak{\ifmmode\let\next\frak@\else
 \def\next{\errmessage{Use \string\frak\space only in math mode}}\fi\next}
\def\goth{\ifmmode\let\next\frak@\else
 \def\next{\errmessage{Use \string\goth\space only in math mode}}\fi\next}
\def\frak@#1{{\frak@@{#1}}}
\def\frak@@#1{\fam\euffam#1}

\def\Bbb{\ifmmode\let\next\Bbb@\else
 \def\next{\errmessage{Use \string\Bbb\space only in math mode}}\fi\next}
\def\Bbb@#1{{\Bbb@@{#1}}}
\def\Bbb@@#1{\fam\msyfam#1}

\def\hexnumber@#1{\ifcase#1 0\or1\or2\or3\or4\or5\or6\or7\or8\or9\or
        A\or B\or C\or D\or E\or F\fi }

\edef\msx@{\hexnumber@\msxfam}
\edef\msy@{\hexnumber@\msyfam}

\mathchardef\blacksquare="0\msx@04
\mathchardef\upharpoonright="3\msx@16

\catcode`\@=12

\begin{document}
\begin{center}
{\Huge Measurable Rectangles}\\
\end{center}

\begin{flushright}
Arnold W. Miller\footnote{Research partially supported by NSF grant
DMS-9024788.}\\
Department of Mathematics\\
University of Wisconsin\\
Madison, WI 53706\\
miller@math.wisc.edu\\
Oct 92, revised Nov 92\\
\end{flushright}

\begin{center}
Abstract.
\end{center}
\begin{quote}
We give an example of a measurable set $E\subseteq\reals$ such that the
set $E^\prime=\{(x,y): x+y\in E\}$ is not in the $\sigma$-algebra
generated by the rectangles with measurable sides.  We also prove
a stronger result that there exists an analytic ($\Sigma^1_1$)
set $E$ such that $E^\prime$ is not
in the $\sigma$-algebra generated by rectangles whose horizontal
side is measurable and vertical side is arbitrary.  The same
results are true when measurable is replaced with property of Baire.
\end{quote}

 The $\sigma$-algebra generated a family ${\cal F}$ of
 subsets of a set $X$ is he smallest family containing
$\cal F$ and closed under taking complements and countable unions.  In
Rao \cite{rao} it is shown that assuming the continuum hypothesis
every subset of the plane $\reals^2$ is in the $\sigma$-algebra
generated by the abstract rectangles, i.e. sets of the form
$A\cross B$ where $A$ and $B$ are arbitrary sets of reals.
In Kunen \cite{kunen} it is shown that it is relatively consistent
with ZFC that not every subset of the plane is in the $\sigma$-algebra
generated by the abstract rectangles.   He shows that this is true
in the Cohen real model.   It also follows from a result of
Rothberger \cite{roth} that if for example $2^{\aleph_0}=\aleph_2$ and
$2^{\aleph_1}=\aleph_{\omega_2}$, then not
not every subset of the plane is in the $\sigma$-algebra
generated by the abstract rectangles.
For a proof of these results see
Miller \cite{millerproj} remark 4 and 5 page 180).

A set is analytic or $\Sigma^1_1$ iff it is the projection of a Borel set.
Answering a question of Ulam, Mansfield \cite{mans1}\cite{mans2} showed
that not every analytic subset of the plane is in the
$\sigma$-algebra generated by the analytic rectangles.   Note that
a rectangle $A\cross B\subseteq \reals \cross \reals$ is analytic iff
both $A$ and $B$ are analytic.

He did this by showing that, in fact, any universal analytic set
is not in the $\sigma$-algebra generated by the rectangles with
measurable sides.  This does the trick because analytic sets are
measurable (see Kuratowski \cite{kur}).  This theorem was also
proved by Rao \cite{rao2}.
Their argument shows a
little more so we give it next.   A set $U\subseteq \reals^2$ is
a universal analytic set
iff it is analytic  and for every analytic set $A\subseteq\reals$ there
exist a real $x$ such that
$$A=U_x=\{y:\langle x,y\rangle\in U\}.$$

\begin{theorem} \label{mans}
  (Mansfield \cite{mans2} and Rao \cite{rao2}) Suppose $U$ is a
  universal analytic set,
  then $U$ is neither in the $\sigma$-algebra generated by rectangles of
  the form
$A\cross B$ with $A\subseteq\reals$ arbitrary and $B\subseteq\reals$
measurable;  nor in the $\sigma$-algebra generated by rectangles of
  the form
$A\cross B$ with $A\subseteq\reals$ arbitrary and $B\subseteq\reals$ having
 the property of Baire.
\end{theorem}
\proof
For any set $U$ in the $\sigma$-algebra generated by rectangles of
  the form
$A\cross B$ with $A\subseteq\reals$ arbitrary and $B\subseteq\reals$
measurable there
is a countable family $\{A_n\cross B_n:n\in\omega\}$ such
that each $B_n$ is measurable and $U$ is in the
$\sigma$-algebra generated by $\{A_n\cross B_n:n\in\omega\}$.

Let $Z$ be a measure zero set and $C_n$ be Borel sets such that
for every $n$ we have $B_n\Delta C_n \subseteq Z$ where $\Delta$ is
the symmetric difference. Since $Z$ is a measure zero set its complement
must contain a perfect set $P$, i.e. a set homeomorphic to the Cantor
space $\cantorspace$. Now for any real $x$ and any set
$V$ in the
$\sigma$-algebra generated by $\{A_n\cross B_n:n\in\omega\}$
we have that $V_x\intersect P$
is Borel.  This is proved by noting that it is trivial if $V=A_n\cross B_n$
(since we have $P\intersect B_n=P\intersect C_n$),
and it is preserved when taking complements and countable unions.
But the set $P$ being perfect must contain a subset $A$ which is analytic
but not Borel.  Since $U$ is universal for some real $x$ we have
that $U_x=A$ and $U_x\intersect P$ is not Borel.  A similar proof
works for the $\sigma$-algebra generated by sets of the form
$A\cross B$ where $B$ has the property of Baire, since if
$B$ has the property of Baire, then for some $G$, an open set,
$B\Delta G$ is meager.
\qed
In Miller \cite{millerrect} it is shown that it is relatively consistent
with ZFC that no universal analytic set is in the
$\sigma$-algebra generated by the abstract rectangles, answering a question
raised by Mansfield.

James Kuelbs
raised the following question\footnote{I want to thank
Walter Rudin for telling me
about this question and also for encouraging me to write
up the solution.}:  If $E\subseteq \reals$ is
measurable, then is
$$E^\prime=\{(x,y) : x+y\in E\}$$
in the $\sigma$-algebra
generated by the rectangles with measurable sides?

One can think of $E^\prime$ as a parallelogram tipped $45$ degrees, so it is
clear by rotation and dilation that $E^\prime$ is measurable if $E$ is.
Note also that since $E^\prime$ is the continuous preimage of $E$,
if $E$ is Borel
then $E^\prime$ is Borel also.

We give a negative answer to Kuelbs' question by showing:

\begin{theorem} For any set $E\subseteq\reals$ we have that
  $E$ is Borel  iff
  $E^\prime$ is in the $\sigma$-algebra generated by rectangles
  of the form $A\cross B$ where  $A$ and $B$ are measurable.
  Similarly,
  $E$ is Borel  iff
  $E^\prime$ is in the $\sigma$-algebra generated by rectangles
  of the form $A\cross B$ where $A$ and $B$ have the property of Baire.
\end{theorem}
\proof
Suppose $E^\prime=\sigma\langle X_n\cross Y_n:n\in\omega\rangle$,
where $X_n$ and $Y_n$ are measurable for every $n$.
Here $\sigma$ is a recipe which describes how a particular set
is built up (using countable intersections and complementation),
i.e., it is the Borel code of $E^\prime$.
Since $X_n$ and $Y_n$ are measurable we can obtain
$A_n$ and $B_n$ Borel sets and $Z$ a Borel set of measure zero
such that:

$X_n\Delta A_n\subseteq Z$  and  $Y_n\Delta B_n\subseteq Z$
 for each $n\in\omega$.

\bigskip

\noindent Claim.  $u\in E$ iff  $\exists x,y\notin Z \; (x+y=u$ and
$\langle x,y\rangle\in \sigma\langle A_n\cross B_n:n\in\omega \rangle)$.

\bigskip

The implication $\leftarrow$ is clear because if $x,y\notin Z$, then
since [$x\in X_n$ iff $x\in A_n$] and [$y\in Y_n$ iff $y\in B_n$] we
have that $(x,y)\in \sigma\langle A_n\cross B_n:n\in\omega \rangle$
iff $(x,y)\in \sigma\langle X_n\cross Y_n:n\in\omega \rangle$
and hence $u\in E$.

The implication $\rightarrow$ is true because of the following
argument. Suppose $u\in E$ is given. Choose $x\notin Z\union (u-Z)$.
Since these are measure zero sets this is easy to do. But now
let $y=u-x$, then $y\notin Z$ since $x\notin u-Z$.  Since
$u\in E$ it must be that $\langle x,y\rangle \in E^\prime$ and
so $\langle x,y\rangle \in  \sigma\langle X_n\cross Y_n:n\in\omega\rangle$
and thus
$\langle x,y\rangle \in  \sigma\langle A_n\cross B_n:n\in\omega\rangle$.
This proves the Claim.

\bigskip
By the Claim,  $E$ is the projection of a Borel set and
hence analytic.  But note that $(E^c)^\prime=(E^\prime)^c$ where
$E^c$ denotes the complement of $E$. It follows that
$E^c$ is also analytic and so by the classical theorem of
Souslin, $E$ is Borel (see Kuratowski \cite{kur}).

A similar proof works for the Baire property since we need only
that every set with the property of Baire differs from some
Borel set by a meager set.

\qed

This answers Kuelbs' question since if $E$ is analytic and not Borel
(or any measurable set which is not Borel), then  $E^\prime$ is
not in the $\sigma$-algebra generated by rectangles with measurable sides.
The argument also shows, for example, that a set $E\subseteq\reals$ is
analytic iff it can be obtained by applying operation A to the
$\sigma$-algebra generated by the rectangles with measurable sides.

\bigskip

Next we show that a slightly stronger result
holds for the sets of the form $E^\prime$.  The following
lemma is the key.

\begin{lemma} \label{mainlem}
  There exists an analytic set $E\subseteq\reals$ such that
  for any $Z$ which
  has measure zero or is meager there exists $x\in\reals$ such
  that $E\setminus(x+Z)$ is not Borel.
\end{lemma}

\proof
Note that we may construct two
such sets, one for category and one for measure, and then putting them
into disjoint  intervals and taking the union would suffice to prove
the lemma.

We first give the proof for category.  We may assume that the set $Z$ is
the countable union of compact nowhere dense sets. This is because,
if $Z^\prime$ is any meager set, then there is such a $Z\supseteq Z^\prime$.
But now if $E\setminus(x+Z)$ is not Borel, then neither is
$E\setminus(x+Z^\prime)$ since
$E\setminus(x+Z)=(E\setminus(x+Z^\prime))\setminus(x+Z)$.

It is a classical result
that the set of irrationals is homeomorphic
to the Baire space, $\bairespace$, which is the space of infinite
sequences of integers with the product topology
(see Kuratowski \cite{kur}).  Let
$h:\reals\setminus\rationals\rightarrow\bairespace$ be a homeomorphism.
For $f,g\in\bairespace$ define $f\leq^* g$ iff for all but finitely
many $n\in\omega$ we have $f(n)\leq g(n)$.  It is not hard to see
the for any countable union of compact sets $F\subseteq \bairespace$
there exists $f\in\bairespace$ such that
$F\subseteq\{g\in\bairespace: g\leq^* f\}$.
Therefore, if $G\subseteq\reals$ is a  $G_\delta$ set
(countable intersection of open sets)
which contains the rationals, then $\reals\setminus G$ is a $\sigma$-compact
subset of $\reals\setminus\rationals$ and therefore there exists
$f\in\bairespace$ such that for every $g\geq^*f$ we have
$h^{-1}(g)\in G$.  (This is a trick going back at least
to Rothberger \cite{roth}).

For $p:\bairespace\arrow\cantorspace$ (the parity function)
by
     $$p(g)(n)\left\{
     \begin{array}{ll} 0 & \mbox{if $g(n)$ is even}\\
                       1 & \mbox{if $g(n)$ is odd}
     \end{array}\right.$$
Let $E_0 \subseteq \cantorspace$ be an analytic set which is not Borel.
Define
 $$E=\{g\in\reals\setminus\rationals:\; p(h(g))\in E_0\}.$$
Now suppose that
$Z\subseteq\reals$ is a meager set which is the countable union
of compact sets and let $G=\reals\setminus Z$.
Let $x_0\in \Intersect_{q\in\rationals}(q-G)$ be arbitrary (this
set is nonempty since each $q-G$ is comeager),
and note that $\rationals\subseteq (x_0+G)$.

Hence there exists $f\in\bairespace$ such
that for every $g\in\bairespace$ with $g\geq^*f$ we have $h^{-1}(g)\in x_0+G$.
Without loss we may assume that for every $n\in\omega$ that $f(n)$ is even.
Let $$P=\{g\in\bairespace:\forall n \; g(n)=f(n)\mbox{ or }g(n)=f(n+1)\}.$$
Then $P$ is homeomorphic to $\cantorspace$ and $Q=h^{-1}(P)\subseteq x+G$
and so $Q\intersect (x+Z)=\emptyset$.  But clearly $Q\intersect E$ is
homeomorphic to $E_0$ via $p\circ h$ and
therefore $E\setminus (x+Z)$ cannot be Borel.

\bigskip

Next we give the proof for measure.
Here we use a coding technique due to Bartoszynski and Judah
\cite{bartjud} who used it to show that a dominating real followed
by a random real gives a perfect set of random reals
(Theorem 2.7 \cite{bartjud}).  We begin by giving the proof in a
slightly different situation, namely $\cantorspace$ instead of
the reals and where $+$ denotes pointwise addition modulo 2
on $\cantorspace$
and the usual product measure on $\cantorspace$.  Afterwards
we will indicate how to modify the
proof for the reals with ordinary addition and Lebesgue measure.

Define $I\subseteq\cantorspace$ to be the
set of all  $x\in\cantorspace$ which have infinitely many ones and infinitely
many zeros.  It is easy to see that $I$ is a $G_\delta$ set.
Any $x\in I$ can be regarded as a sequence
of blocks of consecutive ones, i.e., blocks of consecutive ones each
separated by blocks of one or more zeros.
Define $q:I\arrow\bairespace$ by
$q(x)(n)$ is the length of the $n^{th}$ block of consecutive ones.
As above let $E_0\subseteq\cantorspace$ be an analytic set which
is not Borel and let $p(f)$ for $f\in\bairespace$ be
the parity function.
Let
$$E=\{x\in I : p(q(x))\in E_0\}.$$

We claim that this works, i.e.
given a measure zero set $Z\subseteq\cantorspace$ there exists a real
$x$ such that $E\setminus (x+Z)$ is not Borel.  Let $d_n$ for $n\in\omega$
be the dominating sequence as given in the proof of Theorem 2.7
\cite{bartjud}.  Without loss we may assume that for every
$n\in\omega$ that $d_{2n+2}-d_{2n+1}$ is even if $n$ is even
and odd if $n$ is odd.

According to Bartoszynski and Judah there exists sufficiently random
reals $r,r^\prime\in\cantorspace$ such that following holds.
Let $r^{\prime\prime}$ be defined by
$r^{\prime\prime}(n)=r^{\prime}(n)+1$ mod $2$ for each $n$.  Define
$P$ to be the set of all  $x\in\cantorspace$
such that for every $n$ we have that
$$x\res[d_{2n},d_{2n+1})=r\res[d_{2n},d_{2n+1})$$
and either
$$x\res [d_{2n+1},d_{2n+2})=r^\prime\res[d_{2n+1},d_{2n+2})$$
or
$$x\res [d_{2n+1},d_{2n+2})=r^{\prime\prime}\res[d_{2n+1},d_{2n+2}).$$
The main difficulty of Bartoszynski and Judah's proof is to show that
the perfect set $P$ is disjoint from $Z$. Let $x_1\in\cantorspace$ be defined
by
$$x_1\res[d_{2n},d_{2n+1})=r\res[d_{2n},d_{2n+1})$$
and
$$x_1\res [d_{2n+1},d_{2n+2})=r^\prime\res[d_{2n+1},d_{2n+2}).$$
Let $Q\subseteq\cantorspace$ be the perfect set of all $x\in\cantorspace$ such
that  $x\res [d_{2n},d_{2n+1})$ is constantly zero and
$x\res [d_{2n+1},d_{2n+2})$ is constantly zero or constantly one.
It then follows that $P=x_1+Q$ and that $Q$ is disjoint from $x_1+Z$.
But $E\intersect Q$ is not Borel, because
$d_{2n+2}-d_{2n+1}$ is even if $n$ is even
and odd if $n$ is odd and so it easy to see the $E_0$ is coded into
$E\intersect Q$.

\bigskip

Now we indicate how to modify the above
proof so as to work for the reals with ordinary addition
and Lebesgue measure.  First we modify it to work for the
unit interval $[0,1]$ with ordinary addition modulo one.
Let $s:\cantorspace\rightarrow [0,1]$ be the map defined by
$$s(x)=\Sigma_{n=0}^{\infty}{x(n)\over{2^{n+1}}}.$$
This map is
continuous, measure preserving, and one-to-one except on countably
many points where it is two-to-one.  On the points $x$ where it is
two-to-one let us agree that $s^{-1}(x)$  denotes the preimage of $x$
which is eventually zero.
The main difficulty is that addition mod 1 in $[0,1]$ is
quite different than point-wise addition modulo 2 in $\cantorspace$.
Define for $x,y\in\cantorspace$ the operation $x\oplus y$ to be
$s^{-1}(s(x)+s(y))$ where $s(x)+s(y)$ is the ordinary addition
in $[0,1]$ modulo 1. The operation $\oplus$ just corresponds to
a kind of pointwise addition with carry. Instead of $r^{\prime\prime}$
being the
complement of $r^\prime$ as in the proof of Bartoszynski and Judah
we will take a sparser translate.
Let $Q\subseteq\cantorspace$ be the set of $x\in\cantorspace$ such
that $x(m)=1$ only if for some $n$ we have $m=d_{2n+2}-1$,
i.e., the last element of the interval $[d_{2n+1},d_{2n+2})$.

Note that the set of all
$r\in\cantorspace$ such that for all but finitely many $n$ there exists
$i\in [d_n,d_{n+1})$ such that $r(i)=0$ has measure one.  Hence, by
changing our $d_n$ if necessary we may assume that our random
real $r$ has the property that
for every $n$ there is an $i\in [d_{2n+1},d_{2n+2})$ such that
$r(i)=0$.  This means that when we calculate
$r\oplus x$ for any $x\in Q$ the carry digit on each interval
$[d_{2n+1},d_{2n+2})$ does not propagate out of that interval.
Now by the argument of Bartoszynski and Judah there exists $r\in\cantorspace$
sufficiently random so that  $r\oplus Q$ is disjoint from $s^{-1}(Z)$.

We also use a different coding scheme.   We may assume that for
every $n$ that $d_{2n+2}-1$ is even for even $n$ and odd for odd $n$.
Let $J$ be the set of all $x\in\cantorspace$ such that there are
infinitely many $n$ with $x(n)=1$.  Define $q:\cantorspace\arrow\cantorspace$
be defined by $q(x)=y$ where $\{i_n:n\in\omega\}$ lists in
order all $i$ such that $x(i)=1$ and $y(n)=1$ iff $i_n$ is even.
Now let
$$E=\{x\in [0,1]: q(s^{-1}(x))\in E_0\}.$$
Let $x=s(r)$ and note that $(x+E)\setminus Z$ is not Borel, since
$x+s(Q)$ is disjoint from $Z$. Also note that we can assume
$r(0)=0$ and so $x+s(Q)$ is the same whether we do addition
or addition modulo one. But $E\setminus (-x+Z)$ is just
the translate of $(x+E)\setminus Z$ via $-x$ and so we are done.

\qed

\begin{theorem} \label{main}
There exists $E\subseteq\reals$ which is analytic (hence measurable
and having the property of Baire) such that
$E^\prime=\{(x,y):x+y\in E\}$ is not in either the $\sigma$-algebra
generated by rectangles of the form $A\cross B$
with $A$ arbitrary and $B$ measureable, nor is it in the
$\sigma$-algebra generated by rectangles of the form $A\cross B$
with $A$ arbitrary and $B$ having the property of Baire.
\end{theorem}
\proof

Suppose for contradiction that
$E^\prime=\sigma\langle A_n\cross B_n:n\in\omega \rangle$ where
the $A_n$ are arbitrary and the $B_n$ are measurable.   Let
$Z$ be a measure zero Borel set and $\hat{B_n}$ be Borel such that
$B_n\Delta \hat{B_n}\subseteq Z$ for every $n\in\omega$.
Suppose that $E\setminus (x+Z)$ is not Borel. By translating this
set by $-x$ we must have that $(-x+E)\setminus Z$ is not Borel.
Define $\tilde{B_n}$ as follows.  If $-x\in A_n$ let
$\tilde{B_n}=\hat{B_n}$ and if $-x\notin A_n$ let
$\tilde{B_n}=\emptyset$.  Define
$C=\sigma\langle \tilde{B_n}:n\in\omega\rangle$, i.e., $C$ has
exactly the same Borel code as
$E^\prime=\sigma\langle A_n\cross B_n:n\in\omega \rangle$ except
at the base we substitute $\tilde{B_n}$ for $A_n\cross B_n$.
Since each $\tilde{B_n}$ is Borel the set $C$ is a Borel set.
Now for any $y\notin Z$ we have that
\begin{enumerate}
\item $y\in (-x+E)$ iff
\item $x+y\in E$ iff
\item $(x,y)\in E^\prime$ iff
\item $(x,y)\in \sigma\langle A_n\cross B_n:n\in\omega \rangle$ iff \label{4}
\item $(x,y)\in \sigma\langle A_n\cross \hat{B_n}:n\in\omega\rangle$ iff
\label{5}
\item $y\in \sigma\langle \tilde{B_n}:n\in\omega \rangle=C$. \label{6}
\end{enumerate}
(\ref{4}) and (\ref{5}) are equivalent because $y\notin Z$.
(\ref{5}) and (\ref{6}) are proved equivalent by an easy induction on
the Borel code $\sigma$.

Consequently, for every $y\notin Z$ we have that
$y\in (-x+E)$ iff $y\in C$. But this means that
$(x+E)\setminus Z=C\setminus Z$ which contradicts the assumption
that $(-x+E)\setminus Z$ is not Borel (both $C$ and $Z$ are Borel).
A similar proof works for the property of Baire case.
\qed

There is other work on measurable rectangles which does not seem to
be directly related to this.  For example, Eggleston \cite{egg}
proves that every subset of the plane of positive measure contains
a rectangle $X\cross Y$ with $X$ uncountable (in fact perfect)
and $Y$ of positive
measure.  Martin \cite{mar} gives a metamathematical proof of this
result.

Erdos and Stone \cite{erdst} show that
there exists Borel sets $A$ and $B$ such that the set
$A+B=\{x+y:x\in A, y\in B\}$ is not Borel.

Friedman and Shelah (see Burke \cite{bur} or Steprans \cite{step})
proved that in the Cohen real model for any $F_\sigma$ subset $E$
of the plane, if $E$ contains a rectangle of positive outer measure,
then $E$ contains a rectangle of positive measure.
One corollary of this is that it is consistent that there is
a subset of the plane of full measure which does not contain
any rectangle $A\cross B$ with both $A$ and $B$ having positive
outer measure.  To see this let $E\subseteq\reals$ be  any meager
$F_\sigma$ set with full measure.  Then
$E^\prime=\{(x,y) : x+y \in E\}$ is a subset of the plane of full measure
which is meager.  It cannot contain a rectangle of positive measure
$A\cross B$ since by the classical theorem of Steinhaus the set
$A+B$ would contain an interval and hence $E$ would not be meager.

\end{document}